\documentstyle[12pt,A4]{article}
\author{
	{\bf  Anatol NOWICKI\thanks{Supported by KBN grant 2PO3B130.12;\,
 e-mail: anowicki@omega.wsp.zgora.pl}}\\
	{\it  Institute of Physics, Pedagogical University}\\
	{\it  Plac S{\l}owia\'{n}ski 6, 65-029 Zielona G\'{o}ra, POLAND}\\
	}

\title{Kappa-Deformed Phase Space and Uncertainty Relations}
\date{March 1998}
\newcommand{\ben}{\begin{enumerate}}
\newcommand{\een}{\end{enumerate}}
\newcommand{\be}{\begin{equation}}
\newcommand{\ee}{\end{equation}}
\newcommand{\bea}{\begin{eqnarray}}
\newcommand{\eea}{\end{eqnarray}}
\newcommand{\bc}{\begin{center}}
\newcommand{\ec}{\end{center}}

\begin{document}
\maketitle

\begin{abstract}

We consider two realizations of the $\kappa$-deformed phase space
obtained as a cross product algebra extension of $k$-Poincar\'{e}
algebra. Two kinds of the $\kappa$-deformed uncertainty relations
are briefly discussed.

\end{abstract}

\section{Introduction}
Deformations of space-time symmetries are extensively investigated
in last years. In this approach the notion of symmetries is
generalized  to quantum
groups i.e. Hopf algebras [1]. In this way the Poincar\'{e} symmetry
can be extended to deformed one. Further, we shall consider
the special deformation of space-time symmetry associated to
$k$-deformed Poincar\'{e} algebra with a Hopf algebra structure [2].
In this framework, using the concept of duality, for a given momentum
space one can define configuration space which turns out to be
noncommuting space-time [3]. It appears that both momentum and
configuration spaces form a pair of dual Hopf algebras.\\
It is known that from such  pair of dual algebras one can construct
three different double algebras: Drinfeld double and Drinfeld
codouble, both with Hopf algebra structure, and Heisenberg double
(or more generally, cross product algebra), which is not a Hopf
algebra. Only cross product algebra can be regarded as a
generalization of the standard quantum mechanical phase space,
because in this case we obtain in nondeformed limit the usual
commutation relations between momentum and position operators.\\
Therefore,  deformed phase space cannot be a Hopf algebra.

In chapter 2 we recall the construction of a {\it cross product
algebra} [1], which in our case it is equivalent to
the notion of a Heisenberg double or smash product algebra.

In chapter 3 we show that deformed phase space which is an
extension of the $\kappa$-deformed Minkowski space, dual
to the momentum sector of $\kappa$-Poincar\'{e} algebra [4]
can be defined ambiguously. Roughly speaking, one can define
the commutators between position and momentum in
different ways. In particular, we show that for the
$\kappa$-Minkowski space we can construct two different phase
spaces using the notion of the cross product algebra. In
this way we obtain two different sets of the position
and momentum commutation relations, which generalize in
$\kappa$-deformed way the usual canonical commutation
relations. One type of commutation relations was
considered in [4,5] and the physical implications were discussed
in [6,7].\\ The second form of the commutation
relations give us the uncertainty relations
\begin{displaymath}
\Delta x\Delta p\geq \hbar + \alpha G (\Delta p)^2
\end{displaymath}
(where\quad G- gravitational coupling,\quad $\alpha$-
constant related with the string tension)
\noindent on the algebraic level as a consequence of
the $\kappa$-deformed phase space.
Recently, this kind of uncertainty relations are investigated
in context of quantum gravity and string theories [8,9].
\section{Cross Product Algebra}
Let ${\cal{P}}$ be an algebra, ${\cal{X}}$ a vector space.
\begin{itemize}

\item A {\bf left action} {\it(representation)} of
${\cal{P}}$ on ${\cal{X}}$  is a linear map

\begin{equation}
\triangleright : {\cal{P}}\otimes {\cal{X}}\rightarrow
{\cal{X}} : p\otimes x\rightarrow p\triangleright x
\end{equation}

\noindent such that

\begin{equation}
(p\tilde{p})\triangleright x = p\triangleright (\tilde{p}
\triangleright x)\quad 1\triangleright x = x
\end{equation}

\item We say that ${\cal{X}}$ is a {\bf left
${\cal{P}}$-module}.

\end{itemize}
\noindent In the case where ${\cal{X}}$ is an algebra and
${\cal{P}}$ a bialgebra ({\it or Hopf algebra})
\begin{itemize}

\item We say that ${\cal{X}}$ is a {\bf left
${\cal{P}}$-module algebra} if

\begin{equation}
p\triangleright (x\tilde{x}) = (p_{(1)}\triangleright x) (p_{(2)}
\triangleright \tilde{y}) \quad p\triangleright 1=\epsilon(p) 1
\end{equation}

\noindent where $\epsilon$ denotes counit and we use the
Sweedler's notation
\begin{displaymath}
\Delta(p) = \sum p_{(1)}\otimes p_{(2)}
\end{displaymath}

\end{itemize}
\noindent Let ${\cal{P}}$ be a Hopf algebra and ${\cal{X}}$ a
left ${\cal{P}}$-module algebra.
\begin{itemize}

\item {\bf Left cross product algebra}({\it smash product})
${\cal{X}} >\!\!\!\triangleleft{\cal{P}}$  is a vector
space ${\cal{X}}\otimes {\cal{P}}$ with product [1]

\begin{equation}
(x\otimes p)(\tilde{x}\otimes \tilde{p})=x(p_{(1)}\triangleright
\tilde{x})\otimes p_{(2)}\tilde{p}\qquad
(left\,\,\, cross\,\,\, product)
\end{equation}

\noindent with unit element $1\otimes 1$, where $x, \tilde{x}\in
{\cal{X}}$ and $p, \tilde{p}\in {\cal{P}}$.

\item{\bf commutation relations in a cross product algebra}.
\noindent The obvious isomorphism ${\cal{X}}\sim {\cal{X}}\otimes 1$,
${\cal{P}}\sim
1\otimes {\cal{P}}$ gives us the following cross relations between
the algebras ${\cal{X}}$ and ${\cal{P}}$

\begin{equation}
[x,p]=x\circ p-p\circ x \qquad where \quad x\circ p=x\otimes p
\quad p\circ x=(p_{(1)}\triangleright x)\otimes p_{(2)}
\end{equation}

\end{itemize}
\section{Kappa-Deformed Phase Space as Cross Product
Algebra}
We consider two cases depending on the choice of the momentum sector
basis which give us two different commutation relations between
position and momentum and in consequence yield two kinds of
the uncertainty relations.
\subsection{$\kappa$-Poincar\'{e} Algebra in the Bicrossproduct
Basis}
The $\kappa$-deformed Poincar\'{e} algebra in the bicrossproduct
basis [10] is given by ($g_{\mu\nu}=(-1,1,1,1)$)\\
\noindent - nondeformed ({\it classical}) Lorentz algebra

\begin{equation}
[M_{\mu\nu}, M_{\rho\tau}] = i(g_{\mu\rho}M_{\nu\tau}+
g_{\nu\tau}M_{\mu\rho}-g_{\mu\tau}M_{\nu\rho}
-g_{\nu\rho}M_{\mu\tau})
\end{equation}

\noindent -  deformed covariance relations ($M_i=
\frac12\epsilon_{ijk} M_{jk}, N_i=M_{i0}$)

\begin{equation}
\begin{array}{lll}
[M_i, P_j] & = & i\epsilon_{ijk}P_k\,\qquad [M_i,P_0]\,=\,0\\

[N_i,P_j] & = & i\delta_{ij}\left[\kappa c\sinh(\frac {P_0}{\kappa c})
e^{-\frac{P_0}{\kappa c}}  +
\frac1{2\kappa c} (\vec P)^2\right] - \frac {i}{\kappa c} P_i P_j\\

[N_i,P_0] & = & iP_i\,\qquad [P_\mu, P_\nu]\,=\,0
\end{array}
\end{equation}

\noindent where \quad $\kappa$ - masive deformation parameter and
\quad c - light velocity
and $\kappa$-deformed mass Casimir is given by

\begin{equation}
C_2 =  (2\kappa\sinh \frac{P_0}{2\kappa c} )^2 - {1\over c^2}\vec
P^2 e^{\frac{P_0}{\kappa c}}  = M^2
\end{equation}

\noindent defining coproduct $\Delta$, antipode $S$ and counit
$\epsilon$ as follows

\begin{equation}
\begin{array}{lll}
\Delta(M_i) & = & M_i\otimes 1 + 1\otimes M_i\\

\Delta(N_i) & = & N_i\otimes 1 + e^{-\frac{P_0} {\kappa c}}\otimes
N_i + {\frac{1}{\kappa c}}\epsilon_{ijk}P_j\otimes M_k\\

\Delta(P_0) & = & P_0\otimes 1 + 1\otimes P_0\\

\Delta(P_i) & = & P_i\otimes 1 + e^{-{\frac{P_0}{\kappa c}}}\otimes
P_i \end{array}
\end{equation}

\begin{equation}
S(M_i)=-M_i\qquad S(N_i)=-e^{\frac{P_0}{\kappa c}}N_i+{\frac{1}
{\kappa c}}\epsilon_{ijk}e^{\frac{P_0}{\kappa c}}P_j M_k
\end{equation}

\begin{equation}
S(P_0)=-P_0\qquad S(P_i)=-P_i e^{\frac{P_0}{\kappa c}}
\end{equation}

\begin{equation}
\epsilon(X)=0\qquad for \qquad X= M_i, N_i, P_\mu
\end{equation}

\noindent the $\kappa$-Poincar\'{e} algebra  becomes a Hopf algebra.
It is easy to see that the $\kappa$-Poincar\'{e} algebra contains the
following $\kappa$-deformed Hopf algebra of fourmomentum
${\cal{P}}_\kappa$

\begin{equation}
\begin{array}{lll}
[P_{\mu}, P_{\nu}] & = & 0 \\
\Delta(P_0) & = & P_0 \otimes 1 + 1\otimes P_0 \\
\Delta(P_k) & = & P_k \otimes 1 + e^{-{{P_0}\over \kappa c}}\otimes P_k
\end{array}
\end{equation}

\begin{equation}
S(P_0) = -P_0\qquad S(P_i)=-P_i e^{\frac{P_0}{\kappa c}}\qquad
\epsilon(P_\mu)=0
\end{equation}

\noindent Using the duality relations with second fundamental constant
$\hbar$ (Planck's constant)

\begin{equation}
<x_\mu, P_\nu> = -i\hbar g_{\mu\nu}\qquad g_{\mu\nu} = (-1,1,1,1)
\end{equation}

\noindent we obtain the noncommutative $\kappa$-deformed configuration
space ${\cal{X}}_{\kappa}$ as a Hopf algebra with the following
algebra and coalgebra structure

\begin{equation}
\begin{array}{llll}
[x_{0}, x_{k}] &=& -{{i\hbar}\over \kappa c}x_k\,,\qquad\qquad
& [x_{k}, x_{l}] = 0\\
\Delta(x_\mu ) &=& x_{\mu}\otimes 1 + 1\otimes x_{\mu} &{}\\
S(x_{\mu}) &=&	-x_{\mu}\qquad\qquad
& \epsilon(x_{\mu}) = 0
\end{array}
\end{equation}

\noindent We consider $\kappa$-deformed phase space $\Pi_\kappa \sim
{\cal{X}}_\kappa\otimes {\cal{P}}_\kappa$ as a left cross
product algebra with a product

\begin{equation}
(x\otimes p)(\tilde{x}\otimes \tilde{p})=x(p_{(1)}\triangleright
\tilde{x})\otimes p_{(2)}\tilde{p}
\end{equation}

\noindent and left action

\begin{equation}
p\triangleright x=<p,x_{(2)}>x_{(1)}
\end{equation}

\noindent using the isomorphism $x\sim x\otimes 1$, $p\sim 1\otimes P$
we obtain the commutation relations for $\Pi_\kappa$ in the form

\begin{equation}
\begin{array}{llll}
[x_{0}, x_{k}] &=& -{{i\hbar}\over \kappa c}x_k\,\qquad\qquad
& [x_{k}, x_{l}] = 0\\

[p_\mu, p_\nu] &=& 0 &{}\\

[x_k,p_l] &=& i\hbar \delta_{kl}\qquad &[x_k,p_0]=0\\

[x_0,p_k] &=& {{i\hbar}\over \kappa c} p_k
\qquad &[x_0,p_0]=-i\hbar
\end{array}
\end{equation}

\noindent Introducing the dispersion of the observable
$a$ in quantum mechanical sense by

\begin{equation}
\Delta(a) \ = \ \sqrt{<a^2> - <a>^2} \qquad
\Delta(a)\Delta(b)\geq {1\over 2}|<c>|
\end{equation}

\noindent where \qquad c=[a,b], we obtain $\kappa$-deformed uncertainty
relations in cross product basis $(x_0=ct, E=cp_0)$ [4,5,6]

\begin{equation}
\begin{array}{lll}
\Delta(t)\Delta(x_k) &\geq &{\hbar\over 2\kappa
c^2}|<x_k>|=\frac12\frac{l_\kappa}{c}|<x_k>|\\
\Delta(p_k) \Delta(x_l) &\geq &{1\over 2}\hbar\delta_{kl}\\
\Delta(E) \Delta(t) &\geq & {1\over 2}\hbar\\
\Delta(p_k) \Delta(t) &\geq & {\hbar\over
2\kappa c^2}|<p_k>|=\frac12\frac{l_\kappa}{c} |<p_k>|
\end{array}
\end{equation}

\noindent where $l_\kappa=\frac{\hbar}{\kappa c}$ describes the
fundamental length at which the time variable should already be
considered noncommutative. In the recent estimates $\kappa > 10^{12}
GeV$ therefore	$l_\kappa
< 10^{-26}$cm; in particular one can put $\kappa$ equal to the
Planck mass what implies that $l_\kappa=l_p \simeq 10^{-33}$cm
[see discussion in [6] ).

\subsection{$\kappa$-Poincar\'{e} Algebra in the Standard Basis}
In the standard basis of $\kappa$-Poinar\'{e} algebra [2] the
commutation relations are given by\\
- nondeformed ({\it classical}) $O(3)$-rotation algebra
($M_i=\frac12 \epsilon_{ijk} M_{jk}$)
\begin{equation}
[M_{ij}, M_{kl}] = i(\delta_{ik}M_{jl}+ \delta_{jl}M_{ik}-
\delta_{il}M_{jk} -\delta_{jk}M_{il})
\end{equation}
-  deformed covariance relations
\begin{equation}
\begin{array}{ll}
[M_i, N_j]\,=\, i\epsilon_{ijk}N_k \qquad & [P_\mu,P_\nu]\,=\,0\\

[M_i, P_j]\,=\, i\epsilon_{ijk}P_k \qquad & [M_i,P_0]\,=\,0\\

[N_i, N_j]\,=\,-i\epsilon_{ijk}\left(
M_k \cosh(\frac{P_0}{\kappa c})-{\frac{1}{4(\kappa c)^2}}P_k P_l
M_l\right)\\

[N_i,P_j]\,=\, i\delta_{ij}\kappa c\sinh(\frac {P_0}{\kappa c})
\qquad & [N_i,P_0]\,=\, iP_i
\end{array}
\end{equation}
with $\kappa$-deformed mass Casimir
\begin{equation}
C_2 =  (2\kappa\sinh \frac{P_0}{2\kappa c} )^2 - {1\over c^2}\vec
P^2= M^2
\end{equation}
In this basis the $\kappa$-Poincar\'{e} algebra becomes
a Hopf	algebra if we define coproduct $\Delta$, antipode $S$ and
counit $\epsilon$ in the following way
\begin{equation}
\begin{array}{ll}
\Delta(M_i) \,=& M_i\otimes 1 + 1\otimes M_i\\

\Delta(N_i) \,=& N_i\otimes e^{\frac{P_0}{2\kappa c}} +
e^{-\frac{P_0}{2\kappa c}}\otimes N_i +\\

{} & {\frac{1}{2\kappa c}}\epsilon_{ijk}\left(P_j\otimes
M_ke^{\frac{P_0}{2\kappa c}} + e^{-\frac{P_0}{2\kappa c}}
M_j\otimes P_k\right)\\

\Delta(P_0) \,=& P_0\otimes 1 + 1\otimes P_0\\

\Delta(P_i) \,=& P_i\otimes  e^{-\frac{P_0}{2\kappa c}}  +
e^{\frac{P_0}{2\kappa c}}\otimes P_i
\end{array}
\end{equation}
\begin{equation}
S(M_i)=-M_i\quad S(N_i)=-N_i+{\frac{3i}{2\kappa c}}P_i\quad
S(P_\mu)=-P_\mu
\end{equation}
\begin{equation}
\epsilon(X)=0\qquad for \qquad X= M_i, N_i, P_\mu
\end{equation}
As in the previous case we see, that $\kappa$-Poincar\'{e}
algebra contains the Hopf subalgebra of $\kappa$-deformed
fourmomentum ${\cal{P}}_\kappa$
\begin{equation}
\begin{array}{lll}
[P_{\mu}, P_{\nu}] &=& 0 \\

\Delta(P_0) &=&  P_0 \otimes 1 + 1\otimes P_0\\

\Delta(P_i) &=& P_i\otimes  e^{-\frac{P_0}{2\kappa c}}	+
e^{\frac{P_0}{2\kappa c}}\otimes P_i
\end{array}
\end{equation}
\begin{equation}
S(P_0) = -P_0\qquad S(P_i)=-P_i \qquad \epsilon(P_\mu)=0
\end{equation}
Now, applying the duality relations
we obtain the noncommutative $\kappa$-deformed configuration space
${\cal{X}}_{\kappa}$ as a Hopf algebra with the following algebra
and coalgebra structure
\begin{equation}
\begin{array}{llll}
[x_{0}, x_{k}] &=& -{{i\hbar}\over \kappa c}x_k\,,\qquad\qquad
& [x_{k}, x_{l}]\,=\, 0\\

\Delta(x_\mu ) &=& x_{\mu}\otimes 1 + 1\otimes x_{\mu} &{}\\

S(x_{\mu}) &=&	-x_{\mu}\qquad\qquad
& \epsilon(x_{\mu})\,=\, 0
\end{array}
\end{equation}
It is worthwhile to stress that for both realizations of
bases of fourmomentum sector i.e. bicrossproduct and standard,
we obtain the same commutation relations for
configuration space ${\cal{X}}_\kappa$.
Both bases are related by nonlinear transformation
\begin{displaymath}
P_i\rightarrow P_i e^{\frac{P_0}{2\kappa c}}
\end{displaymath}
Using the left cross product algebra construction we get
commutation relations for $\kappa$-deformed phase space $\Pi_\kappa$
in the standard basis (see also [5])
\begin{equation}
\begin{array}{llll}
[x_{0}, x_{k}] &=& -{{i\hbar}\over \kappa c}x_k\,,\qquad\qquad
& [x_{k}, x_{l}]\,=\, 0\\

[p_\mu, p_\nu] &=& 0 &{}\\

[x_k,p_l] &=& i\hbar \delta_{kl}e^{\frac{P_0}{2\kappa c}}
\qquad &[x_k,p_0]\,=\,0\\

[x_0,p_k] &=& {{i\hbar}\over 2\kappa c} p_k
\qquad &[x_0,p_0]\,=\,-i\hbar
\end{array}
\end{equation}
Therefore, the $\kappa$-deformed uncertainty
relations read $(x_0=ct, E=cp_0)$
\begin{equation}
\begin{array}{lll}
\Delta(t)\Delta(x_k) &\geq &{\hbar\over 2\kappa
c^2}|<x_k>|=\frac12\frac{l_\kappa}{c}|<x_k>|\\

\Delta(p_k) \Delta(x_l) &\geq &{1\over 2}\hbar\delta_{kl}
|<e^{\frac{P_0}{2\kappa c}}>|\\

\Delta(E) \Delta(t) &\geq &{1\over 2}\hbar\\

\Delta(p_k) \Delta(t) &\geq &{\hbar\over
4\kappa c^2}|<p_k>|=\frac14\frac{l_\kappa}{c} |<p_k>|
\end{array}
\end{equation}
and the mass-shell condition
\begin{equation}
(2\kappa\sinh \frac{P_0}{2\kappa c} )^2 - {1\over c^2}\vec P^2= M^2
\end{equation}
This relation one can rewrite as follows
\begin{equation}
exp(\frac{P_0}{2\kappa c})\,=\,\sqrt{\frac{\vec{P}^2 + M^2}
{4\kappa^2 c^2}} + \sqrt{1+\frac{\vec{P}^2 + M^2}{4\kappa^2 c^2}}
\end{equation}
This formula allows us to consider the uncertainty
relations (32) between momentum and position in
the nonrelativistic limit $\quad c\to\infty$
\begin{equation}
\Delta(p_k)\Delta(x_l)>{\hbar\over
4}\delta_{kl}\left[1+\left(1+\frac{M}{2\kappa}\right)^2\right]\,>\,
\frac{\hbar} {2}\delta_{kl}\left(1+\frac{M}{2\kappa}\right)
\end{equation}
Therefore, we see that it is mass dependent relation,
however for $\kappa\gg M$ we get the standard quantum mechanical
uncertainty relation.\\ On the other hand, neglecting the first
term on r.h.s in (34) we have an estimation
\begin{equation}
\Delta(p_k)\Delta(x_l)\,>\,{\hbar\over
2}\delta_{kl}\,\left<\left(1+\frac{\vec{P}^2
+ M^2c^2}{4\kappa^2c^2}\right)^{\frac12}\right>
\end{equation}
using the relation\quad $<\vec{P}^2>\,=\,<\vec{P}>^2 +
(\Delta p)^2$ in the regime \quad $<\vec{P}>^2 + M^2 c^2\ll\kappa^2
c^2,\quad \Delta p\leq\kappa c$ one gets
\begin{equation}
\Delta(p_k)\Delta(x_l)\,>\,{\hbar\over 2}\delta_{kl}\left(1+
\frac{(\Delta p)^2}{8\kappa^2 c^2}\right)
\end{equation}
a modified uncertainty relation which follows from the
analysis of string collisions at Planckian energies (see [4,5]).

\section{Final Remarks}
We considered two realizations of the fourmomentum sector of
$\kappa$-Poincar\'{e} algebra - the {\it bicrossproduct} and
{\it standard} bases. In both cases we discussed the
$\kappa$-deformed phase space obtained by the cross product
algebra construction ({\it Heisenberg double}) and related
uncertainty relations. We see that noncommutativity of space-time
coordinates is the same for both realizations of the phase
space, however the commutation relations between position
and momenta operators are different which give us
two kinds of the uncertainty relations.\\
In the case of the standard basis we showed that the
$\kappa$-deformed uncertainty relations in appropriate limit
give us the modified uncertainty relations postulated in string
theory, now derived on the algebraic level.\\
It is worthwhile to mention here, that the choice of the bases
in the momentum space physically means the choice of
generators which one interpretes as describing physical
momenta. It appears that the choice of the physical momenta
changes radically the uncertainty relations.\bigskip

\begin{center}
{\Large{\bf References}}\bigskip
\end{center}

\begin{enumerate}
\item S. Majid, {\it Foundations of Quantum Group Theory}, (Cambridge
University Press, Cambridge 1995).
\item J. Lukierski, A. Nowicki, H. Ruegg, {\it Phys. Lett.} {\bf B
293} (1993) 419.
\item S. Zakrzewski, {\it Journ. Phys.} {\bf A 27} (1994) 2075;
P. Kosinski, P. Maslanka, {\it The Duality Between
$\kappa$-Poincar\'{e} Algebra and $\kappa$-Poincar\'{e} Group},
(Lodz preprint IM UL 3/1994).
\item J. Lukierski, A. Nowicki, {\it Heisenberg Double Description of
$\kappa$-Poincar\'{e} Algebra and $\kappa$-deformed Phase Space}, in:
{\it Proceedings of XXI International Colloquium on Group Theoretical
Methods in Physics}, (Goslar 1996),\\ q-alg/9702003.
\item A. Nowicki, {\it $\kappa$-Deformed Space-Time Uncertainty
Relations}, in: {\it Proceedings of IX Max Born Symposium},
(Karpacz 1996), q-alg/9702004.
\item G. Amelino-Camelia, J. Lukierski, A. Nowicki,
{\it On the Uncertainty Relations of $\kappa$-Deformed Quantum
Phase Space}, (preprint OUTP-97-40P).
\item G. Amelino-Camelia, {\it Phys. Lett.} {\bf B392} (1997) 283.
\item M. Maggiore, {\it Phys. Lett.} {\bf B319} (1993) 83.
\item A. Kempf, {\it String/Quantum Gravity motivated Uncertainty
Relations\\ and Regularisation in Field Theory},
hep-th/9612082;

A. Kempf, {\it Minimal Lenght Uncertainty Relation and Ultraviolet\\
Regularisation}, preprint DAMTP/96-102, hep-th/9612084.
\item S. Majid, H. Ruegg, {\it Phys. Lett.} {\bf B334} (1994) 348.
\end{enumerate}

\end{document}